\documentclass[11pt,letterpaper]{article}

\usepackage[T2A,T1]{fontenc}
\usepackage[utf8]{inputenc}
\usepackage[russian,english]{babel}

\usepackage{amsmath,amssymb,amsthm}
\usepackage{mathtools}
\numberwithin{equation}{section}

\usepackage{graphicx}
\usepackage{svg}
\usepackage{tikz}
\usetikzlibrary{arrows.meta}

\usepackage{setspace}
\usepackage{xcolor}
\usepackage{multicol}
\usepackage[colorinlistoftodos]{todonotes}

\usepackage{geometry}
\usepackage[numbers]{natbib}

\setcitestyle{numbers,square}
\usepackage[
    unicode,
    colorlinks=true,
    linkcolor=blue,
    urlcolor=blue,
    citecolor=red
]{hyperref}

\newtheoremstyle{mydefstyle}
  {\topsep}          
  {\topsep}          
  {\normalfont}      
  {0pt}              
  {\bfseries}        
  {.}                
  {.5em}             
  {}                 

\theoremstyle{mydefstyle}

\newtheorem{lemma}{Lemma}[section]
\newtheorem{theorem}{Theorem}[section]
\newtheorem{corollary}{Corollary}[section]
\newtheorem{definition}{Definition}[section]
\newtheorem{proposition}{Proposition}[section]

\newtheorem{example}{Example}[section]

\newtheorem{remark}{Remark}[section]



\title{On Pairwise Independent and Strongly Correlated Random Variables}

\author{
    D. V. Gvozdev\textsuperscript{1,2},  A. V. Kolesnikov\textsuperscript{1,2}
    \\[1em]
    \textsuperscript{1}HSE University\\
    \textsuperscript{2}Vega Institute Foundation
}

\date{}

\begin{document}

\bibpunct{[}{]}{,}{n}{,}{,}
\onehalfspacing

\maketitle

\begin{abstract}
Characterizing joint distributions that maximize the centered third-order moment under prescribed marginal distributions and pairwise independence has various statistical applications. This problem can be formulated as a multistochastic optimal transport problem. In this work, we generalize the results on the $(3,2)$-multistochastic Monge--Kantorovich problem to cost functions of the form $c(x,y,z)=f(x)f(y)f(z)$, where $f$ is defined on $[0,1]$ through dyadic expansions:
$f(x)=\sum_{i\geq 1} a_i x_i$, with $x_i\in\{0,1\}$. Under summability and domination conditions on the coefficients, we show that an optimal transport plan is concentrated on the Sierpiński tetrahedron. We also characterize mixtures of geometric coefficient sequences through Hausdorff moment sequences and construct a recursive dual solution with its integral representation for costs with geometric coefficients.
\end{abstract}


\section{Introduction}
It is easy to construct three random variables that are pairwise independent but exhibit nontrivial third-order dependence, for example through elementary constructions involving dice. A more challenging problem is to characterize the joint distribution that maximizes the centered third-order moment when the marginal distributions are given and pairwise independence is imposed. 

This problem can be formulated as a multimarginal optimal transport (MMOT) problem with fixed two-dimensional marginals. For example, consider random variables $X,Y,Z$ uniformly distributed on $[0,1]$ and pairwise independent. A measure $\pi \in \mathcal P_{\lambda}$ can be viewed as a joint distribution of the triple $(X,Y,Z)$ whose two-dimensional marginals are all equal to $\lambda \otimes \lambda$, where $\lambda$ denotes the uniform probability measure on $[0,1]$. We seek a joint distribution $\pi$ that maximizes the centered third-order moment 
$$
\mathbb E_{\pi}\bigl[
(X-\mathbb E[X])(Y-\mathbb E[Y])(Z-\mathbb E[Z])
\bigr].
$$
subject to the prescribed two-dimensional marginals. Uniformity and pairwise independence imply that
$$
\mathbb E_{\pi}\bigl[
(X-\tfrac12)(Y-\tfrac12)(Z-\tfrac12)
\bigr]
=
\mathbb E_{\pi}[XYZ]-\frac18.
$$
Hence maximizing $\mathbb E_{\pi}[XYZ]$ is equivalent to maximizing the centered third-order moment and it can be seen as a multimarginal optimal transportation problem
$$
\sup_{\pi \in \mathcal P_{\lambda}}
    \int_{[0,1]^3} xyz\,d\pi(x,y,z),
$$
This problem was studied in \cite{Gladkov2019} and \cite{GladkovKolesnikovZimin2022}, where it was shown that the support of the optimal joint distribution is exactly the Sierpiński tetrahedron.

Pass \cite{Pass2015} provides a summary of the theory and applications of multimarginal optimal transport. Recent advances include work on joint mixability \cite{KoikeLinWang2024,WangWang2016}, where the authors derive criteria for several classes of distributions under which prescribed marginals admit a coupling satisfying certain moment conditions. A dynamical formulation of MMOT was also recently proposed in \cite{PassShenfeld2025}. From a statistical perspective, one application of MMOT is hypothesis testing: merging several $p$-values for a single hypothesis into one without imposing additional assumptions on their dependence can be formulated as an optimal transport problem (see \cite{VovkWangWang2022,VovkWang2020}). Higher-order dependence also plays a role in statistical testing: the authors of \cite{SejdinovicGrettonBergsma2013} develop nonparametric kernel tests for three-variable interactions. 

Related optimization problems over joint distributions with prescribed marginals arise in dependence modeling and risk aggregation (see, for example, \cite{BernardJiangWang2014,EnnajiMerigotNennaPass2024,PuccettiWang2015}). More specifically, bounds on higher-order product moments and coskewness under dependence uncertainty are investigated in \cite{BernardChenRuschendorfVanduffel2023}.  Optimization over joint distributions with prescribed marginals also appears in econometrics, particularly in partial identification, where optimal transport methods have been used to characterize identified sets in moment models with incomplete data \cite{FanParkPassShi2025} and to derive bounds on causal quantities involving multiple potential outcomes \cite{GaoGeQian2025}.

One possible direction for extending the results in \cite{Gladkov2019, GladkovKolesnikovZimin2022} is to consider other one-dimensional distributions. Quantile transformations allow the corresponding optimization problems to be formulated over the same admissible class $\mathcal P_{\lambda}$ by modifying the cost function. Specifically, for marginal distribution functions $F_X,F_Y,F_Z$, the corresponding cost is
$$
c(u,v,w)=F_X^{-1}(u)F_Y^{-1}(v)F_Z^{-1}(w),
$$
where the inverses are understood as generalized quantile functions. Thus, the prescribed two-dimensional marginals remain uniform product measures, while the choice of one-dimensional distributions is incorporated into the cost.

This observation motivates the following open problem.
\medskip

\noindent \textbf{Problem.} For a bounded nondecreasing function
$f\colon[0,1]\to\mathbb R$, determine the value of
$$
    \sup_{\pi\in\mathcal P_\lambda}
    \int_{[0,1]^3}
    f(x)f(y)f(z)\,d\pi(x,y,z)
$$
and characterize the maximizing plans.
In particular, characterize those functions $f$ for which
the Sierpiński plan remains optimal.\medskip

In this paper, we address this problem for a class of
nondecreasing functions defined through dyadic expansions. We consider cost functions of the form $c(x,y,z)=f(x)f(y)f(z)$, where $f$ is defined on $[0,1]$ through a dyadic expansion:
$
f(x)=\sum_{i\geq 1} a_i x_i,
$
with $x=\sum_{i\geq 1}2^{-i}x_i$ and $x_i\in\{0,1\}$. This class includes the original cost $c(x,y,z)=xyz$ as the special case $a_i=2^{-i}$. Under summability and domination conditions for such functions, it is shown that the optimal solution is also concentrated on the Sierpiński tetrahedron. In addition, we provide a characterization of functions based on geometric sequences through Hausdorff moment sequences, and construct a recursive solution to the dual problem and its integral representation for this class of functions.

The paper is organized as follows. In Section~\ref{section:prelim}, we introduce the preliminary notions needed for the problem. In Section~\ref{section:binary_expansion}, we introduce binary expansion functions and prove their basic properties and study geometric binary expansions and their relation to Hausdorff moment sequences.  In Section~\ref{section:MK-problem}, we study the primal and dual problems. Subsection~\ref{section:MK_primal} proves the existence of an optimal plan supported on the Sierpiński tetrahedron. Subsection~\ref{section:MK_dual} develops the dual formulation for geometric
binary costs, first obtaining a recursive dual potential
and then deriving its integral representation through a
dyadic kernel.

\section{Preliminaries} \label{section:prelim}
Let $X$ and $Y$ be measurable spaces endowed with probability measures $\mu$ and $\nu$, respectively, and let
$c\colon X \times Y \to \mathbb{R}$ be a cost function. The classical \textit{Kantorovich problem}, also known as the transportation problem, consists in minimizing the functional
\begin{equation}
    \inf_{\pi \in \Pi(\mu,\nu)}
    \int_{X \times Y} c(x,y)\,d\pi(x,y),
\end{equation}
where $\Pi(\mu,\nu)$ denotes the set of all probability measures on $X \times Y$ whose marginals are $\mu$ and $\nu$. Equivalently,
$$
    (\operatorname{pr}_X)_\# \pi = \mu,
    \qquad
    (\operatorname{pr}_Y)_\# \pi = \nu.
$$
The elements of $\Pi(\mu,\nu)$ are called \textit{couplings}, or \textit{transport plans}, between $\mu$ and $\nu$. This is the classical Kantorovich formulation of the optimal transportation problem. 
For surveys and comprehensive treatments of the theory, see, for example, 
\cite{AmbrosioBrueSemola2024}, \cite{BogachevKolesnikov2012} or \cite{Villani2009}.

In this paper, we consider an extension of the classical Kantorovich problem,
known as the \textit{multistochastic Monge--Kantorovich problem}. Specifically,
we focus on the $(3,2)$-multistochastic Monge--Kantorovich problem, defined as follows.

\begin{definition}
Let $X_1,X_2,X_3$ be measurable spaces, and let $\mu_{12}$, $\mu_{13}$,
and $\mu_{23}$ be prescribed probability measures on $X_1\times X_2$,
$X_1\times X_3$, and $X_2\times X_3$, respectively. Set
$
    X := X_1\times X_2\times X_3,
$
and denote by $\operatorname{pr}_{ij}\colon X\to X_i\times X_j$
the canonical projection onto coordinates $i$ and $j$. Define
$$
    \mathcal P_{\mu}
    :=
    \left\{
        P\in\mathcal P(X):
        (\operatorname{pr}_{ij})_{\#}P=\mu_{ij}
        \text{ for } 1\leq i<j\leq 3
    \right\},
$$
and assume that $\mathcal P_{\mu}\neq\varnothing$.

Given a measurable cost function
$
    c\colon X\to\mathbb R_{+}\cup\{+\infty\},
$
the \textit{$(3,2)$-multistochastic Kantorovich minimization problem}
consists in finding a measure $P\in\mathcal P_{\mu}$ (called an \textit{optimal transport plan}) such that
$$
    \int_X c\,dP
    =
    \inf_{Q\in\mathcal P_{\mu}}\int_X c\,dQ.
$$
Similarly, the \textit{$(3,2)$-multistochastic Kantorovich maximization problem}
consists in finding a measure $P\in\mathcal P_{\mu}$ such that
$$
    \int_X c\,dP
    =
    \sup_{Q\in\mathcal P_{\mu}}\int_X c\,dQ.
$$
\end{definition}

We extend the results on the $(3,2)$-Monge--Kantorovich problem obtained in
\cite{Gladkov2019, GladkovKolesnikovZimin2022}. The main results of these papers are a characterization of an
optimal plan for the cost function
$c(x,y,z)=xyz$
on the unit cube $[0,1]^3$, under the assumption that the projections onto the
coordinate hyperplanes are given by the two-dimensional Lebesgue measure
$\lambda^2$.

In the present paper, we extend this analysis to a broader class of cost
functions for which the optimal plan has a similar structure. More precisely,
we study the following maximization problem for a broader class of functions:
\begin{equation}
     \label{main_OT_problem}
    \sup_{\pi \in \mathcal P_{\lambda}}
    \int_{[0,1]^3} f(x)f(y)f(z)\,d\pi(x,y,z),
\end{equation}
where $\mathcal P_{\lambda}$ denotes the set of probability measures on
$[0,1]^3$ whose projections onto the coordinate hyperplanes are equal to
$\lambda^2$.

\begin{figure}[h]
    \centering
    \includegraphics[width=0.9\textwidth]{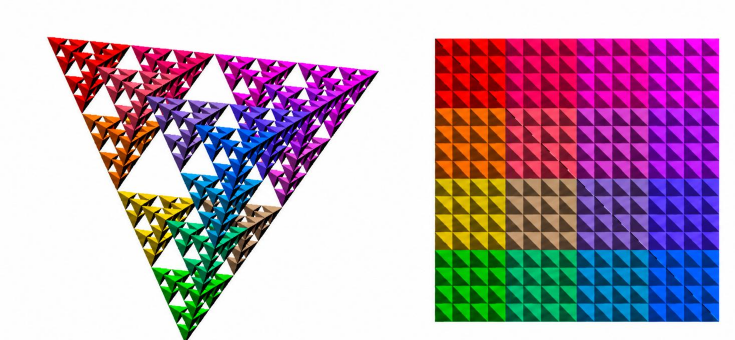}
    \caption{Sierpiński tetrahedron and its two-dimensional projection.}
    \label{fig:sierpinski}
\end{figure}

In particular, we show that the optimal plan is supported on the same
Sierpi\'nski tetrahedron (see Figure~\ref{fig:sierpinski}) and provide a new
formula for a dual solution to the original problem and its extension
to binary expansion costs with geometric coefficients.

\section{Binary Expansion Functions} \label{section:binary_expansion}
\begin{definition}
    \label{def:binary_cost_function}
    Let $(a_i)_{i\geq 1}$ be a sequence of nonnegative numbers satisfying:
\begin{enumerate}
    \item the summability condition
    \begin{equation}
        \tag{Summ}
        \sum_{i=1}^{\infty} a_i < \infty;
        \label{ass:summability}
    \end{equation}

    \item the domination condition
    \begin{equation}
        \tag{Dom}
        a_N > \sum_{n>N} a_n
        \qquad \text{for every } N\geq 1.
        \label{ass:domination}
    \end{equation}
\end{enumerate}
We define a function $f\colon [0,1] \to \mathbb R_+$ by
\begin{equation} 
    f(x):=\sum_{i=1}^{\infty} a_i x_i,
    \label{eq:binary_cost_function}
\end{equation}
where $(x_i)_{i\geq1}$ are the digits in the dyadic (binary) expansion of $x.$
\end{definition}

It turns out that the functions $f$ introduced above have useful structural
properties, such as monotonicity and symmetry with respect to $x=1/2$. The domination
assumption \eqref{ass:domination} implies that the first binary digit at which
two points differ determines the sign of the difference $f(y)-f(x)$.
Moreover, the family of such functions together with zero forms a convex cone.

\begin{proposition} \label{prop:dyadic_cont}
The function
$f$ defined in \eqref{eq:binary_cost_function} is strictly increasing on
$[0,1]$ and continuous $\lambda$-almost everywhere on $[0,1]$.
\end{proposition}

\begin{proof}
We first prove that $f$ is strictly increasing. Let $x,y\in[0,1]$ with
$x<y$. Since the binary expansions of $x$ and $y$ are different, there exists
a first digit at which they differ. Define
$
    N:=\min\{i\geq 1:x_i\neq y_i\}.
$
Since $x<y$, the first differing digit must
satisfy $x_N=0$ and $y_N=1$.
Therefore,
$$
    f(y)-f(x)
    =
    \sum_{i=1}^{\infty}a_i(y_i-x_i)
    =
    a_N+\sum_{i>N}a_i(y_i-x_i).
$$
For every $i>N$, we have $y_i-x_i\geq -1$. 
Consequently, by the domination condition \eqref{ass:domination}, we have
$$
    f(y)-f(x)
    \geq
    a_N-\sum_{i>N}a_i
    >
    0.
$$
Thus, $f(x)<f(y)$ whenever $x<y$. Hence $f$ is strictly increasing on
$[0,1]$.

Since $f$ is monotone, it has at most countably many discontinuities
and is therefore continuous $\lambda$-almost everywhere.
\end{proof}

For our problem it is useful to replace \(f\) by a centered version. Let
$
    A:=\sum_{i=1}^{\infty}a_i.$
Define
\begin{equation}
    F(x):=2f(x)-A
    =\sum_{i=1}^{\infty}a_i(2x_i-1).
    \label{eq:symmetric_F}
\end{equation}
Then $f(x)=\frac{F(x)+A}{2}.$
The function $F$ is nondecreasing because $f$ is nondecreasing. Moreover,
away from dyadic points, the binary expansion of \(1-x\) is obtained by
complementing all binary digits of \(x\). Hence \(F\) is antisymmetric with
respect to \(x=1/2\) outside a countable set.

\begin{proposition}
\label{prop:F_properties}
The function $F$ has the following properties:
\begin{enumerate}
    \item[\textit{(i)}] $F$ is nondecreasing.
    \item[\textit{(ii)}] For every non-dyadic $x\in(0,1)\colon
        F(1-x)=-F(x).
    $
    \item[\textit{(iii)}] Fix \(n\geq 1\) and define the $n$-th tail by
    $$
        F_n(u):=
        \sum_{k=1}^{\infty}a_{n+k-1}(2u_k-1).
    $$
    Then
    \begin{itemize}
        \item $\displaystyle F_n(u)<0
        \quad \text{if } u\in(0,1/2)$;
        \item $\displaystyle F_n(u)>0
        \quad \text{if } u\in(1/2,1).$
    \end{itemize}
\end{enumerate}
\end{proposition}
\begin{proof}
\textit{(i)} The monotonicity of $F$ follows directly from the monotonicity of $f$.

\textit{(ii)} Take any non-dyadic $x\in(0,1)$. Note that $1-x$ has the complementary binary expansion $(1-x)_i=1-x_i$; then
    $$
    F(1-x) = \sum_{i=1}^{\infty}a_i (2(1-x_i)-1)=\sum_{i=1}^{\infty}a_i (1-2x_i) = -F(x).
    $$

    \textit{(iii)} If $u_1=0$, then we have by \eqref{ass:domination} 
    $$
    F_n(u) \leq -a_n + \sum_{k>n} a_k <0.
    $$
    If $u_1=1,$ then 
    $$
    F_n(u) \geq a_n - \sum_{k>n} a_k >0,
    $$
    because $a_n > \sum_{k>n} a_k$ for any $n \geq 1.$
\end{proof}

\subsection{Geometric Binary Expansions}
We now discuss a special case of the binary functions introduced above.
For a given $a\geq2$, define
\begin{equation}
 \label{eq:binary_geomcost_function}
    C_a(x):=\sum_{i=1}^{\infty}\frac{x_i}{a^i},
\end{equation}
where $(x_i)_{i\geq 1}$ are the digits in the dyadic expansion of $x$.
The endpoint $a=2$ gives the linear function $C_2(x)=x$, up to the
choice of binary expansion at dyadic points. For parameters strictly greater than 2, it is easy to check that $C_a$ satisfies \eqref{ass:summability} and
\eqref{ass:domination}.
\begin{definition}
\label{def:binary_geomcost_function}
Define the cone generated by the family $\{C_a\}_{a\geq 2}$, defined in
\eqref{eq:binary_geomcost_function}, by
$$
\mathcal C
:=
\left\{
\sum_{a\in A}\lambda_a C_a \;:\;
A\subset [2,\infty) \text{ is at most countable},\ 
\lambda_a\geq 0 \text{ for all } a\in A,\ 
\sum_{a\in A}\lambda_a<\infty
\right\},
$$
where the sum is understood pointwise.
\end{definition}

\begin{theorem}
\label{theorem:conic_cost_function}
Let $f$ be the uniform limit of a sequence $(f_n)\subset \mathcal C$
such that, writing
$$
f_n=\sum_{a\in A_n}\lambda_{n,a}C_a,
$$
we have
$ 
\sup_n\sum_{a\in A_n}\lambda_{n,a}<\infty.
$
Then $f$ admits a representation
\begin{equation}
    \label{conic_cost_function} 
    f(x)=\sum_{i=1}^{\infty}c_i x_i,
\qquad
c_i:=\int_{[0,1/2]}s^i\,d\mu(s),
\end{equation}
for some finite nonnegative Borel measure $\mu$ on $[0,1/2]$.
\end{theorem}

\begin{proof}
First consider an arbitrary function $f\in\mathcal C$. Then
$$
f(x)
=
\sum_{a\in A}\lambda_a C_a(x)
=
\sum_{a\in A}\lambda_a
\sum_{i=1}^{\infty}\frac{x_i}{a^i},
\qquad
\lambda_a\geq 0,
\qquad
\sum_{a\in A}\lambda_a<\infty .
$$
Since all terms are nonnegative, Fubini-Tonelli's theorem allows us to exchange the
order of summation. Hence
$$
f(x)
=
\sum_{i=1}^{\infty}x_i
\sum_{a\in A}\frac{\lambda_a}{a^i}.
$$
Define the finite atomic measure $\nu$ on $[2,\infty)$ by
$$
\nu:=\sum_{a\in A}\lambda_a\delta_a, \qquad
\sum_{a\in A}\frac{\lambda_a}{a^i}
=
\int_{[2,\infty)} t^{-i}\,d\nu(t).
$$
Let $T(t):=1/t$ and define the push-forward measure $\mu:= \nu \circ T^{-1}$ on $(0,1/2]$. By the change-of-variables formula,
$$
\int_{[2,\infty)} t^{-i}\,d\nu(t)
=
\int_{(0,1/2]} s^i\,d\mu(s)
=:c_i.
$$
Therefore every $f\in\mathcal C$ admits the representation
$$
f(x)=\sum_{i=1}^{\infty} c_i x_i,
\qquad
c_i=\int_{(0,1/2]} s^i\,d\mu(s).
$$

Finally, let $\mu_n$ be the measures constructed above. Their masses are uniformly bounded by assumption. Compactness of $[0,1/2]$ therefore gives a subsequence converging weakly to a finite nonnegative measure $\mu$. For each binary sequence $(x_i)_{i\geq1}$, the function
$\sum_{i=1}^\infty x_i s^i
$
is continuous on $[0,1/2]$, since its series converges uniformly. Hence
$$
f(x)=\lim_{j\to\infty}f_{n_j}(x)
=\lim_{j\to\infty}\int_{(0,1/2]} \sum_{i=1}^\infty x_is^id\mu_{n_j}
=\sum_{i=1}^\infty x_i\int_{(0,1/2]} s^i\,d\mu(s) = \sum_{i=1}^\infty c_ix_i,
$$
as we sought. 
\end{proof}
\begin{corollary} \label{corollary:Hausdorf_MMS}
The sequence $(c_i)_{i\geq 1}$ in \eqref{conic_cost_function} is a
\textit{completely monotone sequence}; that is,
$$
    (-1)^k(\Delta^k c)_n \geq 0
    \qquad \forall n\geq 1,\ k\geq 0,
$$
where $\Delta$ is the forward difference operator defined by
$
    (\Delta c)_n := c_{n+1}-c_n.
$
Moreover, for every $n\geq 1$, we have
$
    c_n \geq 2c_{n+1}.
$
\end{corollary}




\begin{remark}
Since $\mu$ is a finite nonnegative measure on $[0,1/2]$, the sequence
$(c_i)_{i\geq 1}$ can be viewed as a rescaled Hausdorff moment sequence.
Hausdorff moment sequences are precisely the completely monotone sequences;
see \cite{Hausdorff1921}.
\end{remark}

\section{Monge--Kantorovich Problem with Binary Expansion Costs} \label{section:MK-problem}
\subsection{Primal Problem} \label{section:MK_primal}
We consider the following $(3,2)$-Kantorovich maximization problem
\begin{equation}
    \int_{[0,1]^3} f(x)f(y)f(z)\,d\pi
    \longrightarrow
    \max_{\pi\in\mathcal P_\lambda},
    \label{eq:primal_f_problem}
\end{equation}
where $f$ is defined in \eqref{eq:binary_cost_function}.

The use of $\max$, rather than $\sup$, is justified by a compactness argument. The admissible set $\mathcal P_\lambda$ is weakly compact, since it is tight and closed under weak convergence (the proof is similar to the standard Kantorovich problem, see, for example, \cite{BogachevKolesnikov2012} for more details). The cost function $f(x)f(y)f(z)$ is bounded. Its discontinuities are contained in the union of coordinate hyperplanes corresponding to dyadic points. Every admissible plan has one-dimensional Lebesgue marginals; hence every such hyperplane has zero mass under every $\pi\in\mathcal P_\lambda$. Therefore, by the standard weak convergence theorem for bounded functions whose discontinuity set has zero mass, the objective functional is continuous on $\mathcal P_\lambda$. Consequently, the supremum is attained, and we might replace $\sup$ with $\max$ in the problem.

The symmetric normalization allows us to reduce this problem to the maximization of the cubic term involving $F$ defined in \eqref{eq:symmetric_F}.

\begin{lemma}
\label{lem:reduction_to_F}
Problem \eqref{eq:primal_f_problem} is equivalent to
\begin{equation}
    \int_{[0,1]^3} F(x)F(y)F(z)\,d\pi
    \longrightarrow
    \max_{\pi\in\mathcal P_\lambda}.
    \label{eq:primal_F_problem}
\end{equation}
\end{lemma}
\begin{proof}
    Since $f=(F+A)/2$, we have
$$
    f(x)f(y)f(z)
    =
    \frac{1}{8}(F(x)+A)(F(y)+A)(F(z)+A).
$$
Expanding the product gives a sum of constant, linear, quadratic, and cubic
terms in $F$. The integrals of the constant term, the linear terms, and the
quadratic terms are fixed over $\mathcal P_\lambda$, because the
two-dimensional marginals of every admissible plan are fixed and equal to
$\lambda^2$. Therefore, the only term depending on $\pi$ is
$$
    \int_{[0,1]^3}F(x)F(y)F(z)\,d\pi.
$$
Thus, maximizing \eqref{eq:primal_f_problem} is equivalent to maximizing
\eqref{eq:primal_F_problem}.
\end{proof}
\begin{remark}
    The exceptional dyadic points form a countable set. Since every admissible
plan has Lebesgue one-dimensional marginals, every coordinate hyperplane of
the form $\{x=r\}$, $\{y=r\}$, or $\{z=r\}$, with $r$ dyadic, has zero
mass. Thus, the values of $F$ at dyadic points do not affect the value of the
transport problem.
\end{remark}

We can now prove the analog of the support result for the cost
$c(x,y,z)=f(x)f(y)f(z)$.

\begin{theorem} \label{cost_a:theorem}
There exists an optimal plan for \eqref{eq:primal_f_problem} supported on the
Sierpi\'nski tetrahedron.
\end{theorem}
\begin{proof}
By Lemma~\ref{lem:reduction_to_F}, it is enough to consider the symmetric
problem \eqref{eq:primal_F_problem}.

Next, note that the problem is invariant under the transformations
$$
T_{xy}(x,y,z)=(1-x,1-y,z),\qquad
T_{xz}(x,y,z)=(1-x,y,1-z),\qquad
T_{yz}(x,y,z)=(x,1-y,1-z),
$$
each of which pushes forward any $\mu\in\mathcal P_\lambda$ to a measure in $\mathcal P_\lambda$. By symmetry, it suffices to consider $T_{xy}$. Denote $\mu^{xy}:=\mu\circ T_{xy}^{-1}$. Then, by Proposition~\ref{prop:F_properties}(ii),
$$
\int_{[0,1]^3} F(x)F(y)F(z)\, d\mu^{xy}
=\int_{[0,1]^3} F(1-x)F(1-y)F(z)\, d\mu
=\int_{[0,1]^3} F(x)F(y)F(z)\, d\mu.
$$
Thus the objective is invariant under $T_{xy}$ (and similarly under $T_{xz},T_{yz}$). In particular, if $\pi$ solves \eqref{eq:primal_f_problem}, then so do $\pi^{xy},\pi^{xz},\pi^{yz}$, and therefore the symmetrized plan
$$
\pi_1:=\frac{\pi+\pi^{xy}+\pi^{xz}+\pi^{yz}}{4}
$$
is optimal as well.

Now we decompose $[0,1]^3$ into two sets $I_1$ and $I_2$ by
$$
I_2 :=
\left[0,\tfrac12\right]^3
\cup
\left(\left[\tfrac12,1\right]^2 \times \left[0,\tfrac12\right]\right)
\cup
\left(\left[0,\tfrac12\right] \times \left[\tfrac12,1\right]^2\right)
\cup
\left(\left[\tfrac12,1\right] \times \left[0,\tfrac12\right] \times \left[\tfrac12,1\right]\right),
$$
$$
I_1 := \overline{[0,1]^3\setminus I_2}.
$$
Both $I_1$ and $I_2$ are invariant under $T_{xy},T_{xz},T_{yz}$, hence the restricted measures
$$
\pi_{I_1}:=\pi_1|_{I_1},\qquad \pi_{I_2}:=\pi_1|_{I_2}
$$
are invariant as well.

We claim that $\pi_{I_2}=0$. Let $\pi_1 = \pi_{I_1} + \pi_{I_2}.$ Assume that $\pi_{I_2}\neq 0$. Let $T_x(x,y,z) := (1-x,y,z)$ and define the measure
$
\widetilde\pi := \pi_{I_1} + \pi_{I_2}\circ T_x^{-1}.
$
Since projections are symmetric with respect to $x\mapsto 1-x$, the map $T_x$ sends $\mathcal P_\lambda$ to itself.
Moreover, $I_1$ and $I_2$ are unions of dyadic half-cubes, and $T_x(I_2)=I_1$.
Therefore $\widetilde\pi\in\mathcal P_\lambda$.
Since $\widetilde\pi$ and $\pi_1$ have the same unreflected component on $I_1$, we have
\begin{align*}
    \int_{[0,1]^3} F(x)F(y)F(z) d\widetilde\pi -\int_{[0,1]^3} F(x)F(y)F(z) d\pi_{1}
= \int_{I_1} F(x)F(y)F(z) d(\pi_{I_2}\circ T_x^{-1}) \\
   -\int_{I_2} F(x)F(y)F(z) d\pi_{I_2}.
\end{align*}
Using Proposition~\ref{prop:F_properties}(ii), we get
$$
\int_{I_1} F(x)F(y)F(z) d(\pi_{I_2}\circ T_x^{-1})
= \int_{I_2} F(1-x)F(y)F(z) d\pi_{I_2}
= -\int_{I_2} F(x)F(y)F(z) d\pi_{I_2}.
$$
Hence
$$
\int_{[0,1]^3} F(x)F(y)F(z) d\widetilde\pi -\int_{[0,1]^3} F(x)F(y)F(z) d\pi_{1}
= -2\int_{I_2} F(x)F(y)F(z) d\pi_{I_2}.
$$

On $I_2$ we have $F(x)F(y)F(z)<0$ by Proposition~\ref{prop:F_properties}(iii),
except possibly on the boundary where at least one coordinate equals $1/2$ (the value of $F(1/2)$ is immaterial).
Since every $\mu\in\mathcal P_\lambda$ has one-dimensional marginals equal to $\lambda$, we obtain
$$
\pi_{I_2}\Big(\{x=\tfrac12\}\cup\{y=\tfrac12\}\cup\{z=\tfrac12\}\Big)=0,
$$
so the integrand is $<0$ $\pi_{I_2}$-a.e. on $I_2$. Because $\pi_{I_2}\neq 0$, it follows that
$$
\int_{I_2} F(x)F(y)F(z) d\pi_{I_2} < 0,
$$
and therefore $\widetilde \pi$ attains a strictly larger value of the objective, contradicting the optimality of $\pi_1$. Thus, we have $\pi_{I_2}=0$.

The same argument can now be repeated inside each of these four subcubes. After
rescaling a retained subcube to $[0,1]^3$, the relevant part of the function $F$ is described by the shifted tail $F_2$, which satisfies the same sign
properties as $F$ by Proposition~\ref{prop:F_properties}(iii). Therefore one can
construct an optimal plan supported on the union of $4^2$ dyadic subcubes of
side length $2^{-2}$.

Iterating this construction, we obtain a decreasing sequence of compact sets
$
    J_1\supset J_2\supset \cdots,
$
where $J_n$ is a union of $4^n$ dyadic cubes of side length $2^{-n}$, and
an optimal plan $\pi_n$ supported on $J_n$. By compactness of
$\mathcal P_\lambda$, a subsequence of $(\pi_n)$ converges weakly to some
$\pi^\ast\in\mathcal P_\lambda$. Since the cost is bounded and its
discontinuity set has zero mass for every admissible plan, the objective is
continuous along such weakly convergent sequences. Hence $\pi^\ast$ is still
optimal.

Finally,
$$
    \operatorname{supp}\pi^\ast
    \subset
    \bigcap_{n=1}^{\infty}J_n.
$$
The intersection
$
    J:=\bigcap_{n=1}^{\infty}J_n
$
is precisely the Sierpi\'nski tetrahedron. Therefore, there exists an optimal
plan supported on the Sierpi\'nski tetrahedron.
\end{proof}

In analogy with the results of \cite{Gladkov2019}, we may formulate our theorem in terms of a transport map $T$ and a solution to the minimization problem.

\begin{corollary}
    The pushforward measure $\pi := \lambda^2 \circ T^{-1}$ of the two-dimensional Lebesgue measure $\lambda^2$ under the map
    $
        T\colon (x,y) \mapsto \bigl(x,y,1-(x \oplus y)\bigr),
    $
    where $\oplus$ denotes bitwise addition modulo 2, is a solution to problem \eqref{eq:primal_f_problem}.
    If, instead of maximizing the total cost functional, one considers the corresponding minimization problem, then the optimal transport map is given by
    $
        T(x,y) = (x,y,x\oplus y).
    $
\end{corollary}

\begin{corollary} \label{cor:geom_cost_sierpinski}
An optimal plan $\pi$ for
\begin{equation}
\label{main_geom_OT_problem}
    \int_{[0,1]^3} C_a(x)C_a(y)C_a(z)\,d\pi
    \to \max_{\pi \in \mathcal P_\lambda}
\end{equation}
is supported on the Sierpi\'nski tetrahedron.
\end{corollary}

The following corollary follows from Theorem~\ref{cost_a:theorem} and  Corollary~\ref{corollary:Hausdorf_MMS}, since every $f \in \mathcal C$ admitting a representation of the form \eqref{conic_cost_function} with a finite nonnegative measure satisfies \eqref{ass:summability} and \eqref{ass:domination}.

\begin{corollary} \label{theorem:conic_cost_function_plan}
Let $\pi$ be the plan concentrated on the Sierpi\'nski tetrahedron, then $\pi$ is also optimal for
$$
    \int_{[0,1]^3} f(x)f(y)f(z)\, d\pi
    \longrightarrow
    \max_{\pi\in\mathcal P_{\lambda}}
$$
for every $f\in\overline{\mathcal C}$, where the closure is taken in the uniform norm.
\end{corollary}

The following example shows that the sequence $(a_i)_{i\geq1}$ in Definition~\ref{def:binary_cost_function} need not be completely monotone, unlike the coefficient sequences in Corollary~\ref{theorem:conic_cost_function_plan}.

\begin{example}
     Consider the sequence
    $$
    a_1=1,\qquad a_2=0.6,\qquad a_3=0.1,\qquad
    a_i=\frac{10^{-2}}{4^{i-4}} \quad \text{for all } i\geq 4.
    $$
    Then $(a_i)_{i\geq1}$ is summable. Moreover, the domination condition
    \eqref{ass:domination} holds:
    $$
    \sum_{i\geq 4} a_i
    =
    10^{-2}\sum_{j=0}^{\infty}\frac{1}{4^j}
    =\frac{1}{75}
    <
    0.1=a_3,  \quad 2a_3<a_2, \quad a_2+2a_3<a_1,
    $$
    and for every $N\geq 4$ we have
    $
    a_N
    >
    \sum_{i>N} a_i,
    $
    since the tail from $N$ onward is a geometric sequence.

    However, the sequence is not completely monotone. Indeed,
    $$
    (\Delta^2 a)_1
    =
    a_3-2a_2+a_1
    =
    0.1-1.2+1
    =
    -0.1<0,
    $$
    whereas complete monotonicity would require $(\Delta^2 a)_1\geq 0$.
\end{example}

\subsection{Dual Problem} \label{section:MK_dual}
For this subsection, we consider a minimization analogue of problem \eqref{main_OT_problem}, so that the formulas are a bit shorter. Namely, we consider
$$
    \int_{[0,1]^3} f(x)f(y)f(z)\,d\pi
    \longrightarrow \min_{\pi\in\mathcal P_\lambda},
$$
where $\pi$ has Lebesgue projections onto the principal coordinate hyperplanes. We may formulate the corresponding dual problem:
$$
    \int_{[0,1]^2} F(x,y)\,dxdy
    +
    \int_{[0,1]^2} G(y,z)\,dydz
    +
    \int_{[0,1]^2} H(z,x)\,dzdx
    \longrightarrow \max,
$$
subject to the pointwise constraint
$$
    F(x,y)+G(y,z)+H(z,x)
    \leq
    f(x)f(y)f(z)
    \qquad \text{for all } (x,y,z)\in[0,1]^3.
$$
Indeed, since every admissible plan $\pi\in\mathcal P_\lambda$ has fixed
two-dimensional marginals equal to $\lambda^2$, any admissible triple
$(F,G,H)$ gives
$$
\begin{aligned}
    \int_{[0,1]^3} f(x)f(y)f(z)\,d\pi
    &\geq
    \int_{[0,1]^3}
    \bigl[
        F(x,y)+G(y,z)+H(z,x)
    \bigr]\,d\pi  \\
    &=
    \int_{[0,1]^2}F(x,y)\,dxdy
    +
    \int_{[0,1]^2}G(y,z)\,dydz
    +
    \int_{[0,1]^2}H(z,x)\,dzdx.
\end{aligned}
$$
Thus the dual problem provides lower bounds for the primal minimization
problem.

By the symmetry of the cost $f(x)f(y)f(z)$ and of the marginal constraints,
one may reduce the dual problem to the symmetric case
$$
    F=G=H,
    \qquad
    F(x,y)=F(y,x).
$$
In this form the dual problem becomes
$$
    3\int_{[0,1]^2}F(x,y)\,dxdy
    \longrightarrow \max,
$$
subject to
$$
F(x,y)+F(y,z)+F(z,x)
    \leq
    f(x)f(y)f(z)
    \qquad\text{for all }(x,y,z)\in[0,1]^3.
$$

By the standard Kantorovich duality argument, if a symmetric function
$F$ satisfies
\begin{equation}
 \label{problem:duality_ineq}
F(x,y)+F(y,z)+F(z,x)
    \leq
    f(x)f(y)f(z)
    \qquad\text{for all }(x,y,z)\in[0,1]^3,
\end{equation}
and equality holds $\pi^\ast$-almost everywhere, namely
\begin{equation}
 \label{problem:duality_eq}
F(x,y)+F(y,z)+F(z,x)
    =
    f(x)f(y)f(z)
    \qquad \pi^\ast\text{-a.e.},
\end{equation}
then $\pi^\ast$ is optimal for the primal minimization problem.

\subsubsection{Recursive Solution for Cost Functions Generated by Geometric Binary Expansions}
Consider the cost functions defined in \eqref{eq:binary_geomcost_function}. Our goal is to construct a symmetric function $F$ satisfying
\begin{equation} 
    F(x,y)+F(y,z)+F(z,x)
    \leq
    C_a(x)C_a(y)C_a(z),
    \label{eq:lower_dual_ineq}
\end{equation}
for all $x,y,z \in [0,1]$, with equality on the support $z=x\oplus y$. 

Note that, similarly to the primal problem, we can consider, without loss of generality, a centered version of $C_a$ for a given $a\geq2$ defined as $$g(x):= 2C_a(x) - A, \quad A = \sum_{i=1}^\infty a^{-i} = \frac{1}{a-1}.$$

We introduce a shift operator $T$ on binary expansions $\Omega:=\{0,1\}^\mathbb N$ defined as
$$
Tx :=  \begin{cases}
    (x_2,x_3,\cdots), & x_1=1\\
(1-x_2, 1-x_3, \cdots), & x_1=0.
\end{cases}
$$

The following theorem characterizes a possible dual solution $H$ for costs of the form $c(x,y,z)=g(x)g(y)g(z)$ using self-similarity properties of $g$.

\begin{theorem}
    For a given $a\geq 2$, define
$$
\varphi(u,v):= uv + \frac{u+v}{2} + \frac13.
$$
and
\begin{equation}
\label{eq:folded_H_series}
H(x,y):=\sum_{k=1}^{\infty}
\frac{
\varphi \bigl(g(T^k x),g(T^k y)\bigr)
}{a^{3k}}
\end{equation}
Then $-H$ is a symmetric dual solution to the $(3,2)$-multistochastic Monge-Kantorovich minimization problem with the cost $ c(x,y,z)=g(x)g(y)g(z)$.
\end{theorem}
\begin{proof}
Since $|g|\leq A$ and
    $
    |\varphi(u,v)|\leq A^2+A+1/3   
    $ for all $|u|,|v|\leq A$,
    the series \eqref{eq:folded_H_series} converges uniformly.
    Thus $H$ is bounded and symmetric, and
    \begin{equation}
        \label{eq:folded_H_recursion}
        H(x,y)
        =
        \frac{
            H(Tx,Ty)+\varphi(g(Tx),g(Ty))
        }{a^3}.
    \end{equation}

We now prove the dual inequality. Define the bounded function
$$  
D(x,y,z):=
H(x,y)+H(y,z)+H(z,x)-g(x)g(y)g(z).
$$
For brevity, put
$
u:=g(Tx), v:=g(Ty), w:=g(Tz).
$
Since
$$
\varphi(u,v)+\varphi(v,w)+\varphi(w,u)
=
(1+u)(1+v)(1+w)-uvw,
$$
equation \eqref{eq:folded_H_recursion} gives
\begin{align*}
H(x,y)+H(y,z)+H(z,x)
&=
\frac{
D(Tx,Ty,Tz)+(1+u)(1+v)(1+w)
}{a^{3}}\\
&=
\frac{D(Tx,Ty,Tz)}{a^{3}}
+|g(x)g(y)g(z)|,
\end{align*}
where the second equality follows from
the self-similarity of $g$: $|g(x)|=\frac{1+g(Tx)}a.$
Therefore,
\begin{equation}
\label{eq:folded_defect_recursion}
D(x,y,z)
=
\frac{D(Tx,Ty,Tz)}{a^3}
+|g(x)g(y)g(z)|-g(x)g(y)g(z).
\end{equation}

Iterating this identity $N$ times, we obtain
\begin{align*}
D(x,y,z)
=\frac{D(T^Nx,T^Ny,T^Nz)}{a^{3N}} + \sum_{k=0}^{N-1}
\frac{|g(T^kx)g(T^ky)g(T^kz)|
-g(T^kx)g(T^ky)g(T^kz)}{a^{3k}}
\end{align*}
The first term tends to zero because $D$ is bounded.
Every summand is nonnegative, so $D\geq0$.
Using
\eqref{eq:folded_defect_recursion}, we conclude that
$$
H(x,y)+H(y,z)+H(z,x)
\geq |g(x)g(y)g(z)|.
$$
In particular,
$
-H(x,y)-H(y,z)-H(z,x)
\leq g(x)g(y)g(z).
$

Next we establish equality on $z=x\oplus y$.
For $\varepsilon\in\{0,1\}$, define
$$
\mathcal D_\varepsilon
:=
\left\{
(x,y,z)\in\Omega^3:
x_j\oplus y_j\oplus z_j=\varepsilon
\text{ for every }j\geq1
\right\}.
$$
The digit formula for the shift is given by
$
(Tx)_j=x_{j+1}\oplus(1-x_1).
$
We have
$$
(Tx)_j\oplus(Ty)_j\oplus(Tz)_j
=
x_{j+1}\oplus y_{j+1}\oplus z_{j+1}\oplus1\oplus x_1\oplus y_1\oplus z_1 = \varepsilon \oplus 1 \oplus \varepsilon=1.
$$
Thus, the coordinatewise map $(T,T,T)$ sends both
$\mathcal D_0$ and $\mathcal D_1$ into $\mathcal D_1$.

On $\mathcal D_1$ Proposition~\ref{prop:F_properties}(iii) implies $
g(x)g(y)g(z)\geq0
$ since the number of zero first digits is even. Note that $\mathcal D_1$ is invariant under $(T,T,T)$,
every remainder term in the iterated identity vanishes
on this set, and we obtain
$$
D(x,y,z) = \frac{D(T^Nx,T^Ny,T^Nz)}{a^{3N}}.
$$
Letting $N\to\infty$ and using the boundedness of $D$, we conclude that
$D(x,y,z)=0$ on $\mathcal D_1$.

If $(x,y,z)\in\mathcal D_0$, then Proposition~\ref{prop:F_properties}(iii) implies
$
g(x)g(y)g(z)\leq0.
$
Note that, after applying the shift, we have $(Tx,Ty,Tz)\in\mathcal D_1$, and hence
$D(Tx,Ty,Tz)=0$. Consequently,
$$
H(x,y)+H(y,z)+H(z,x)
=
|g(x)g(y)g(z)|
=
-g(x)g(y)g(z).
$$
Since $\mathcal D_0$ is precisely the set $z=x\oplus y$,
this proves equality for the centered cost on the required set.
\end{proof}

Note that the original function $C_a$ can be expressed as $C_a = (g+A)/2$. This yields the following corollary.

\begin{corollary}
    For a given $a\geq 2$, the original dual problem \eqref{eq:lower_dual_ineq} admits the following dual solution $F$:
\begin{equation} \label{eq:recursive_formula}
    F(x,y):=
\frac18\left[
-H(x,y)
+A g(x)g(y)
+\frac{A^2}{2}\bigl(g(x)+g(y)\bigr)
+\frac{A^3}{3}
\right].
\end{equation}
\end{corollary}

\subsubsection{Integral Representation of the Dual Solution}
In \cite{GladkovKolesnikovZimin2022}, it is stated that a dual solution can be represented as an integral with respect to a positive measure on $[0,1]$. In this section, we show that such a representation also holds for our problem.

The main idea behind obtaining the integral representation is to introduce  a finite-increment analogue of taking a mixed derivative $\Delta F(x,y)$ to eliminate all terms of the form $g(x)+g(y)$. Since $C_a$ is right-continuous and increasing, it defines a Lebesgue--Stieltjes measure $\mu_a$ by 
\begin{equation} \label{eq:leb-stiel}
    \mu_a((u,v]):=C_a(v)-C_a(u) \text{ for every $(u,v]\subseteq[0,1]$,} \quad \mu_a(\{0\})=0.
\end{equation}
For $a>2$, this measure is positive and purely atomic. Using the symmetry on dyadic sets, we express $\Delta F$ as an integral of a kernel with respect to $\mu_a\otimes\mu_a$. Finally, we recover the terms depending on only one variable and obtain the integral representation stated in Theorem~\ref{theorem:integral_rep}.

\begin{proposition}
    Let $a>2$, and let $\mu_a$ be the Lebesgue--Stieltjes measure defined in \eqref{eq:leb-stiel}. Its restriction to $(0,1)$ coincides with the distributional
derivative $dC_a$:
    \begin{equation}
  \label{eq:weak_der_Ca}
        \mu_a
        =
        \frac{a-2}{a-1}
        \sum_{n=1}^{\infty}a^{-n}
        \sum_{\substack{1\le m\le 2^n-1\\ m\text{ odd}}}
        \delta_{m/2^n}.
    \end{equation}
\end{proposition}
\begin{proof}
    Since $C_a$ is increasing, its distributional derivative $\mu_a$
    is a positive measure of total mass
    $$
        \mu_a((0,1))
        =C_a(1-)-C_a(0+)=\frac{1}{a-1}.
    $$
    At a dyadic point $r=m/2^n$, with $m$ odd, the two binary
    expansions differ by replacing the $n$-th digit $0$ followed
    by all ones with $1$ followed by all zeros. Thus
    $$
        \mu_a(\{r\})
        =C_a(r+)-C_a(r-)
        =a^{-n}-\sum_{k>n}a^{-k}
        =\frac{a-2}{a-1}a^{-n}.
    $$
    There are $2^{n-1}$ such points at level $n$, and their total mass is
    $$
        \frac{a-2}{a-1}
        \sum_{n=1}^{\infty}2^{n-1}a^{-n}
        =\frac{1}{a-1}
        =\mu_a((0,1)).
    $$
    Hence, these atoms exhaust $\mu_a$, proving the representation \eqref{eq:weak_der_Ca}.
\end{proof}

Now define the rectangular increment $\Delta$ for functions of the form $B(x,y)$ by
\begin{equation} \label{eq:inc_DeltaF}
    \Delta B(x,y) := B(x,y) -B(x,0)-B(0,y)+B(0,0).
\end{equation}
It's easy to check that for a product we have
$
\Delta (h(x)h(y))= (h(x)-h(0))(h(y)-h(0)).
$
By linearity, we get 
$$
\Delta F(x,y) = -\frac{1}{8}\Delta H + \frac{A}{8} \Delta (g(x)g(y)) =-\frac{1}{8}\Delta H + \frac{A}{2} C_a(x)C_a(y),
$$
since $g(x)=2C_a(x)-A$. The same rule applies to function $H$ defined in \eqref{eq:folded_H_series} and we get
\begin{equation} \label{eq:Delta_F}
    \Delta F(x,y) = \frac{A}{2}C_a(x)C_a(y) - \frac{1}{8} \sum_{k=1}^\infty \frac{(h_k(x)-h_k(0))(h_k(y)-h_k(0))}{a^{3k}},
\end{equation}
where we defined $h_k(x):=g(T^kx)$.

The measure $\mu_a$ is supported on the dyadic points,
each of which has a unique representation $s=p/2^n$ with $p$ odd.
Define the reduced dyadic level by
$$
n(s)=
\begin{cases}
    n, & \text{if } s=p/2^n \text{ with } p \text{ odd},\\
    +\infty, & \text{if } s \text{ is non-dyadic}.
\end{cases}
$$
The following observation, which follows from the symmetries of the shift $T$, is the key to obtaining an integral representation.
\begin{lemma}
    Let $a>2$, and put $h_k(x):=g(T^kx)$ for $k\geq1$.
    Then, for every $x\in[0,1]$,
    \begin{equation} \label{eq:hk_integral}
        h_k(x)-h_k(0)
        =2a^k\int_{(0,x]}\rho_k(s)\,d\mu_a(s),
    \end{equation}
    where, for $s\in(0,1)$,
    \begin{equation} \label{eq:rho_k}
        \rho_k(s):=
        \begin{cases}
            2s_k-1, & k<n(s),\\
            0, & k\geq n(s),
        \end{cases} \qquad \rho_k(0)=\rho_k(1)=0.
    \end{equation}
\end{lemma}

\begin{proof}
    The digit formula for the $T$ gives
    $
        (T^kx)_j=x_{k+j}\oplus(1-x_k),
    $
    and hence
    $
        h_k(x)=(2x_k-1)g(L^kx).
    $
    where $L$ denotes the left-shift operator.
    
    On each open dyadic cell of level $k$, the first $k$ digits
    are fixed, and
    $$
        C_a(L^kx)
        =a^k\left(C_a(x)-\sum_{i=1}^k a^{-i}x_i\right).
    $$
    Since $g=2C_a-A$, it follows that
    $
        dh_k=2a^k(2x_k-1)dC_a.
    $ It remains to examine the cell boundaries. At a dyadic point
    of reduced level $m\leq k$, the binary expansions corresponding
    to the left and right limits are
    $$
        b_1\cdots b_{m-1}0\,111\cdots,
        \qquad
        b_1\cdots b_{m-1}1\,000\cdots.
    $$
    If $m=k$, then the formula
    $h_k(x)=(2x_k-1)g(L^kx)$ gives the value $-A$ on both sides,
    because $g(0)=-A$ and $g(1)=A$.
    If $m<k$, it gives the value $A$ on both sides.
    Thus $h_k$ has no jump at any boundary of a $k$-level cell.

    Thus, we get
    $
        dh_k=2a^k\rho_kd\mu_a  
    $ on $(0,1)$ where $\rho_k$ is defined above.  Integrating the derivative measure over $(0,x]$ yields
    $$
        h_k(x)-h_k(0)
        =2a^k\int_{(0,x]}\rho_k(s)\,d\mu_a(s),
    $$
    as we sought.
\end{proof}

Now we are ready state and prove our representation results.

\begin{theorem} \label{theorem:integral_rep}
    For given $a>2$, the dual problem \eqref{eq:lower_dual_ineq} admits the following integral representation 
    $$
    F(x,y) = I_a(x,y) - \frac{1}{4}I_a(x,x) - \frac14 I_a(y,y),
    $$
    where $I_a(x,y)=\Delta F(x,y)$ and can be expressed using a kernel $K_a$ as
    $$
    I_a(x,y) := \int_{(0,x] \times (0,y]} K_a(s,t) d\mu_a(s) d\mu_a(t), \quad K_a(s,t) := \frac{1}{2} \sum_{k=1}^\infty \frac{1-\rho_k(s)\rho_k(t)}{a^k},
    $$
    and $\rho_k$ are defined in \eqref{eq:rho_k} and $\mu_a$ is a Lebesgue--Stieltjes measure from \eqref{eq:leb-stiel}.
\end{theorem}
\begin{proof}
    First, we show that $I_a(x,y)=\Delta F(x,y)$. Note that
    $$
    C_a(x)C_a(y) = \int_{(0,x]\times(0,y]} d\mu_a(s)d\mu_a(t),
    $$
    which, combined with substitution of \eqref{eq:hk_integral} into \eqref{eq:Delta_F} gives
    \begin{align*}
        \Delta F(x,y) &=\frac{A}{2}\int_{(0,x]\times(0,y]} d\mu_a(s)d\mu_a(t) - \frac{1}{2} \sum_{k=1}^\infty a^{-k} \left( \int_{(0,x]}\rho_k(s)\,d\mu_a(s)\right) \cdot \left(\int_{(0,y]}\rho_k(t)\,d\mu_a(t)\right) \\&= \int_{(0,x]\times(0,y]} \left(\frac{A}{2} - \frac{1}{2}\sum_{k=1}^\infty a^{-k} \rho_k(s)\rho_k(t) \right) d\mu_a(s)d\mu_a(t),
    \end{align*}
    where exchanging the series and the integrals is justified by the fact that $a,|\rho_k|\leq 1$ and $\mu_a$ is finite.

    Finally, note that $A = \sum_{k=1}^\infty a^{-k}$ and  
    $$
    \Delta F(x,y) = \int_{(0,x]\times(0,y]} \left(\frac{1}{2}\sum_{k=1}^\infty a^{-k} (1-\rho_k(s)\rho_k(t)) \right) d\mu_a(s)d\mu_a(t)= \int_{(0,x]\times(0,y]} K_a(s,t) d\mu_a(s)d\mu_a(t),
    $$
    as we sought.

    Finally, to recover $F$, note that points $(0,0,0)$ and $(x,x,0)$ belong to a $\pi$-support and have zero-costs. Thus,
    $
    F(0,0)=0$, and $F(x,x) +2F(x,0)=0.
    $
    It follows that 
    $$
    I_a(x,x) = \Delta F(x,x) = F(x,x)-2F(x,0)+F(0,0) = -4F(x,0).
    $$
    Substituting this into the definition of $\Delta F(x,y)$, we obtain
    $$
     F(x,y) = I_a(x,y) - \frac{1}{4}I_a(x,x) - \frac14 I_a(y,y), 
    $$
    as we sought.
\end{proof}
\begin{remark}
    For $a=2$ we have $C_2(x) = x$ and $dC_2 = dx$. For non-dyadic $s,t$, 
    $$
    K_2(s,t) = \sum_{k=1}^\infty 2^{-k}\frac{1-\rho_k(s)\rho_k(t)}{2} = \sum_{k=1}^\infty 2^{-k}(s_k \oplus t_k) = s\oplus t \quad \text{a.e.}
    $$
    Thus, the representation reduces to 
    $$
    I_2(x,y) = \int_{0}^x \int_0^y (s\oplus t) dtds,
    $$
    which is the integral representation for the original problem in \cite{ GladkovKolesnikovZimin2022}.
\end{remark}
\medskip\medskip

\begin{small}
\noindent \textbf{Acknowledgements.} The second author gratefully acknowledges the support of the Basic Research Program of HSE University. 
\end{small}

\medskip\medskip \begin{small}
\noindent \textbf{AI Disclosure Statement.} The dual solutions were developed with the assistance of ChatGPT 6 Astra (OpenAI). We edited and checked all generated proofs for accuracy. Any remaining errors are the responsibility of the authors. 
\end{small}

\bibliographystyle{amsplain}
\bibliography{references}

@article{Gladkov2019,
  author  = {Gladkov, Nikita A. and Kolesnikov, Alexander V. and Zimin, Alexander P.},
  title   = {On Multistochastic {Monge--Kantorovich} Problem, Bitwise Operations, and Fractals},
  journal = {Calculus of Variations and Partial Differential Equations},
  year    = {2019},
  volume  = {58},
  number  = {5},
  pages   = {173},
  doi     = {10.1007/s00526-019-1610-4}
}

@book{Villani2009,
  author    = {Villani, C{\'e}dric},
  title     = {Optimal Transport: Old and New},
  series    = {Grundlehren der mathematischen Wissenschaften},
  volume    = {338},
  publisher = {Springer},
  address   = {Berlin, Heidelberg},
  year      = {2009},
  doi       = {10.1007/978-3-540-71050-9}
}

@article{GladkovKolesnikovZimin2022,
  author  = {Gladkov, Nikita A. and Kolesnikov, Alexander V. and Zimin, Alexander P.},
  title   = {The Multistochastic Monge--Kantorovich Problem},
  journal = {Journal of Mathematical Analysis and Applications},
  year    = {2022},
  volume  = {506},
  number  = {2},
  pages   = {125666},
  doi     = {10.1016/j.jmaa.2021.125666}
}

@article{Hausdorff1921,
  author  = {Hausdorff, Felix},
  title   = {Summationsmethoden und Momentfolgen {I-II}},
  journal = {Mathematische Zeitschrift},
  year    = {1921},
  volume  = {9},
  pages   = {74--109, 280--299}
}

@article{BogachevKolesnikov2012,
  author  = {Bogachev, Vladimir I. and Kolesnikov, Alexander V.},
  title   = {The Monge--Kantorovich Problem: Achievements, Connections, and Perspectives},
  journal = {Russian Mathematical Surveys},
  year    = {2012},
  volume  = {67},
  number  = {5},
  pages   = {785--890},
  doi     = {10.1070/RM2012v067n05ABEH004808}
}

@book{AmbrosioBrueSemola2024,
  author    = {Ambrosio, Luigi and Bru{\'e}, Elia and Semola, Daniele},
  title     = {Lectures on Optimal Transport},
  series    = {UNITEXT},
  volume    = {169},
  edition   = {2},
  publisher = {Springer Cham},
  year      = {2024},
  doi       = {10.1007/978-3-031-76834-7}
}

@article{PuccettiWang2015,
  author  = {Puccetti, Giovanni and Wang, Ruodu},
  title   = {Extremal Dependence Concepts},
  journal = {Statistical Science},
  year    = {2015},
  volume  = {30},
  number  = {4},
  pages   = {485--517},
  doi     = {10.1214/15-STS525}
}

@article{BernardJiangWang2014,
  author  = {Bernard, Carole and Jiang, Xiao and Wang, Ruodu},
  title   = {Risk Aggregation with Dependence Uncertainty},
  journal = {Insurance: Mathematics and Economics},
  year    = {2014},
  volume  = {54},
  pages   = {93--108},
  doi     = {10.1016/j.insmatheco.2013.11.005}
}

@article{BernardChenRuschendorfVanduffel2023,
  author  = {Bernard, Carole and Chen, Jinghui and
             R{\"u}schendorf, Ludger and Vanduffel, Steven},
  title   = {Coskewness under Dependence Uncertainty},
  journal = {Statistics \& Probability Letters},
  year    = {2023},
  volume  = {199},
  pages   = {109853},
  doi     = {10.1016/j.spl.2023.109853}
}

@inproceedings{SejdinovicGrettonBergsma2013,
  author    = {Sejdinovic, Dino and Gretton, Arthur and Bergsma, Wicher},
  title     = {A Kernel Test for Three-Variable Interactions},
  booktitle = {Advances in Neural Information Processing Systems},
  year      = {2013},
  volume    = {26},
  publisher = {Curran Associates, Inc.},
  url       = {https://papers.nips.cc/paper_files/paper/2013/hash/076a0c97d09cf1a0ec3e19c7f2529f2b-Abstract.html}
}

@misc{FanParkPassShi2025,
  author        = {Fan, Yanqin and Park, Hyeonseok and
                   Pass, Brendan and Shi, Xuetao},
  title         = {Partial Identification in Moment Models with
                   Incomplete Data---A Conditional Optimal Transport Approach},
  year          = {2025},
  eprint        = {2503.16098},
  archivePrefix = {arXiv},
  primaryClass  = {econ.EM},
  note          = {Preprint, version 2},
  url           = {https://arxiv.org/abs/2503.16098}
}

@inproceedings{GaoGeQian2025,
  author    = {Gao, Zijun and Ge, Shu and Qian, Jian},
  title     = {Bridging Multiple Worlds: Multi-marginal Optimal Transport
               for Causal Partial-identification Problem},
  booktitle = {Proceedings of the 28th International Conference
               on Artificial Intelligence and Statistics},
  year      = {2025},
  volume    = {258},
  series    = {Proceedings of Machine Learning Research},
  pages     = {721--729},
  publisher = {PMLR},
  url       = {https://proceedings.mlr.press/v258/gao25b.html}
}

@article{Pass2015,
  author  = {Pass, Brendan},
  title   = {Multi-marginal optimal transport: Theory and applications},
  journal = {ESAIM: Mathematical Modelling and Numerical Analysis},
  year    = {2015},
  volume  = {49},
  number  = {6},
  pages   = {1771--1790},
  doi     = {10.1051/m2an/2015020}
}

@article{EnnajiMerigotNennaPass2024,
  author  = {Ennaji, Hamza and M{\'e}rigot, Quentin
             and Nenna, Luca and Pass, Brendan},
  title   = {Robust risk management via multi-marginal optimal transport},
  journal = {Journal of Optimization Theory and Applications},
  year    = {2024},
  volume  = {202},
  number  = {2},
  pages   = {554--581},
  doi     = {10.1007/s10957-024-02438-x}
}

@article{WangWang2016,
  author  = {Wang, Bin and Wang, Ruodu},
  title   = {Joint mixability},
  journal = {Mathematics of Operations Research},
  year    = {2016},
  volume  = {41},
  number  = {3},
  pages   = {808--826},
  doi     = {10.1287/moor.2015.0755}
}

@article{KoikeLinWang2024,
  author  = {Koike, Takaaki and Lin, Liyuan and Wang, Ruodu},
  title   = {Joint mixability and notions of negative dependence},
  journal = {Mathematics of Operations Research},
  year    = {2024},
  volume  = {49},
  number  = {4},
  pages   = {2786--2802},
  doi     = {10.1287/moor.2022.0121}
}

@article{VovkWang2020,
  author  = {Vovk, Vladimir and Wang, Ruodu},
  title   = {Combining {$p$}-values via averaging},
  journal = {Biometrika},
  year    = {2020},
  volume  = {107},
  number  = {4},
  pages   = {791--808},
  doi     = {10.1093/biomet/asaa027}
}

@article{VovkWangWang2022,
  author  = {Vovk, Vladimir and Wang, Bin and Wang, Ruodu},
  title   = {Admissible ways of merging {$p$}-values
             under arbitrary dependence},
  journal = {The Annals of Statistics},
  year    = {2022},
  volume  = {50},
  number  = {1},
  pages   = {351--375},
  doi     = {10.1214/21-AOS2109}
}

@misc{PassShenfeld2025,
  author        = {Pass, Brendan and Shenfeld, Yair},
  title         = {A dynamical formulation of multi-marginal optimal transport},
  year          = {2025},
  eprint        = {2509.22494},
  archivePrefix = {arXiv},
  primaryClass  = {math.OC},
  note          = {arXiv preprint arXiv:2509.22494},
  doi           = {10.48550/arXiv.2509.22494},
  url           = {https://arxiv.org/abs/2509.22494}
}
\end{document}